\documentclass[11pt, reqno, oneside]{amsproc}

\newcommand{\golab}{Go\l \c{a}b}

\usepackage{amsmath, amsthm, amssymb, amsfonts, graphicx, overpic, bm}

\newcommand{\area}[1]{\operatorname{area}(#1)}
\newcommand{\perimeter}[1]{\operatorname{per}(#1)}
\DeclareMathOperator{\bd}{\partial\!}
\DeclareMathOperator{\conv}{conv}

\newcommand{\normal}{\dashv}
\newcommand{\notnormal}{\not\,\dashv}
\newcommand{\antinormal}{\normal_{\mathrm{a}}}

\newcommand{\Ray}[1]{\overrightarrow{\,#1\,\,}}

\newcommand{\epsi}{\varepsilon}
\newcommand{\fhi}{\varphi}

\newcommand{\norm}[1]{\Vert#1\Vert}
\newcommand{\antinorm}[1]{\|#1\|_{\mathrm{a}}}
\newcommand{\mnorm}[1]{\vert#1\vert}

\newcommand{\R}{\mathbf{R}}

\newcommand{\abs}[1]{\vert#1\vert}
\newcommand{\bigabs}[1]{\big\vert#1\big\vert}

\newcommand{\myangle}{\sphericalangle}

\newtheorem{theorem}{Theorem}
\newtheorem{proposition}{Proposition}
\newtheorem{lemma}{Lemma}
\newtheorem{corollary}{Corollary}

\begin{document}

\title{Antinorms and Radon curves}
\thanks{Research supported by a grant from an agreement between the Deutsche Forschungsgemeinschaft in Germany and the National Research Foundation in South Africa. Parts of this paper were written during a visit of the first author to the Department of Mathematical Sciences of the University of South Africa.}
\author{Horst Martini}
\address{Fakult\"at f\"ur Mathematik,
        Technische Universit\"at Chemnitz,
        D-09107 Chemnitz, Germany}
\email{\texttt{martini@mathematik.tu-chemnitz.de}}
\author{Konrad J. Swanepoel}
\address{Department of Mathematical Sciences,
        University of South Africa, PO Box 392,
        Pretoria 0003, South Africa}
\email{\texttt{swanekj@unisa.ac.za}}
\date{}
\maketitle


\section{Introduction}
In this paper we consider two notions that have been discovered and rediscovered by geometers and analysts since 1917 up to the present day.
The first is that of a \emph{Radon curve}, introduced by Radon in 1917 \cite{Radon}.
It is a special kind of centrally symmetric closed convex curve in the plane.
Any centrally symmetric closed convex curve in the plane defines a norm turning the plane into a two-dimensional normed space or \emph{Minkowski plane}.
Finite-dimensional normed spaces or \emph{Minkowski spaces} were introduced by Minkowski in \cite{Minkowski}, and a special case (the $L_p$ norm for $p=4$) was even alluded to in Riemann's famous \emph{Habilitationsvortrag} \cite{Riemann}.
For a general introduction to Minkowski spaces, see Thompson's book \cite{MR97f:52001} and the surveys \cite{MR2002h:46015, MSII}.
The norm thus defined by a Radon curve is called a \emph{Radon norm}, and the corresponding Minkowski plane a \emph{Radon plane}.
Radon planes have many remarkable, almost-Euclidean properties.
For a survey on Radon planes, including further results, see \cite[Section~6]{MR2002h:46015}.

The second notion is that of an \emph{antinorm}.
This is a norm dual in a certain sense to the norm of an arbitrary Minkowski plane.
It is a special case of the \emph{Minkowski content} of a set in a Minkowski space introduced by Minkowski \cite{Minkowski2}.
Busemann \cite{MR9:372h, MR17:779a} showed that the circles in the antinorm (\emph{anticircles}) are the solutions to the isoperimetric problem in a Minkowski plane.
He also showed that anticircles are circles (i.e.\ the antinorm is proportional to the norm) exactly when the circles are Radon curves.
The plane with the antinorm turns out to be isometric to the dual normed space of the plane with the original norm.
Note however that the antinorm and the norm are defined in the same plane.
Since there is no \emph{natural} way of identifying a vector space and its dual, even in the finite-dimensional case, any identification of the dual normed plane plane with the antinorm must involve some arbitrary choice.
In principle one would have to choose an invertible linear transformation, or equivalently, fix a coordinate system (four degrees of freedom in the two-dimensional case).
Choosing a Euclidean structure and using polarity is also sufficient  (three degrees of freedom).
However, we will explain how one only has to choose a unit of area and an orientation (one degree of freedom), since in the plane there exists up to a constant factor only one symplectic bilinear form (a multiple of the determinant) \cite{MR90h:51003}.
It is in this almost-natural context that our constructions will be done.
Since area enters the picture, it is not surprising that the antinorm ties together the norm and the area, such as in the isoperimetric problem mentioned above, as well as in other results.

Many known results in Euclidean geometry also hold for Radon planes, for example the triangle and parallelogram area formulas, certain theorems on angular bisectors, the area formula of a polygon circumscribed about a circle, certain isoperimetric inequalities, and the non-expansive property of certain non-linear projections.
These results may be further generalized to arbitrary Minkowski planes if we formally change the statement of the result by referring in some places to the antinorm instead of the norm.
It is the purpose of this mainly expository paper to give a list of results on antinorms that generalize results true for Radon norms, and in many cases characterize Radon norms among all norms in the plane.
Many of the results are old, well-known, and have often been rediscovered.
However, for most of the results we give streamlined proofs.
Also, some of the characterizations of Radon curves seem not to have appeared previously in print.
In particular we mention that Corollaries~\ref{radonbisector}, \ref{newcor}, and \ref{zenodorus} are new.

\section{Fundamentals}

\subsection{Convex curves}
By a \emph{plane} we mean a two-dimensional real vector space $V$.
Everything in this paper will be a part of affine geometry (see \cite[chapter~13]{Coxeter}), and our presentation will be as coordinate-free as possible.
A \emph{convex body} in $V$ is a convex closed bounded set with interior points.
A \emph{convex curve} is the boundary $\bd C$ of a convex body $C$.
We denote the area of convex body $C$ by $\area{C}$.
We need the following two technical lemmas.
The first is proved e.g.\ in \cite{MR47:5732}.
The second follows from the fact that a convex body equals the intersection of its supporting half spaces.

\begin{lemma}\label{tech1}
Let $C_1$ and $C_2$ be two convex bodies in $V$ both containing the origin in their interiors.
Suppose that for any non-zero vector $v\in V$ we have that the points where the ray $\{\lambda v:\lambda\geq 0\}$ intersects $C_1$ and $C_2$ have parallel supporting lines.
Then $C_1=\lambda C_2$ for some $\lambda>0$.
\end{lemma}

\begin{lemma}\label{tech2}
Let $C_1$ and $C_2$ be two convex bodies in the plane such that each supporting line of $C_1$ is also a supporting line of $C_2$.
Then $C_1=C_2$.
\end{lemma}

The following lemma will help us to deduce results for convex bodies from results for polygons.
It can be proved easily using compactness.
\begin{lemma}\label{tech3}
Given a convex body $C$ with supporting lines $\ell_1,\dots,\ell_n$.
For any $\epsi>0$ there exists a polygon $P$ circumscribed to $C$ whose sides are contained in $\ell_1,\dots,\ell_n$ as (possibly degenerate) sides, and such that $P\subseteq (1+\epsi)C+v$ for some $v\in V$.
\end{lemma}

\subsection{Normed planes and dual planes}

As usual, a \emph{norm} on $V$ is a real-valued function $\norm{\cdot}$ on $V$ satisfying
\begin{itemize}
\item $\norm{x}\geq 0$,
\item $\norm{x}=0\implies x=o$, 
\item $\norm{\lambda x}=\abs{\lambda}\norm{x}$, and 
\item $\norm{x+y}\leq\norm{x}+\norm{y}$.
\end{itemize}
The \emph{unit ball} is \[B=B_{\norm{\cdot}} := \{x\in V:\norm{x}\leq 1\},\] and the \emph{unit circle} \[\bd B := \{x\in V:\norm{x}=1\}.\]
Equivalently, one may start with a bounded, centrally symmetric convex body $B$ with non-empty interior, and define \[\norm{x}=\inf\{\lambda>0 : \lambda^{-1}x\in B\}.\]
The \emph{perimeter} or length of a convex body $C$ in the norm is denoted by $\perimeter{C}$.
It can be defined as the limit of the perimeters (in the norm) of polygons inscribed in $C$.
Note that even though we will introduce a second norm (the antinorm) we will always measure the perimeter of $C$ in the original norm unless the contrary is clearly stated, to avoid any confusion.

A \emph{functional} on $V$ is a linear transformation $\fhi:V\to\R$.
The \emph{dual plane} of $V$ is the set $V^\ast$ of all functionals on $V$.
Then $V^\ast$ is also a two-dimensional real vector space.
The \emph{dual norm} is defined by
\[ \norm{\fhi}^\ast := \sup\{\fhi(x):\norm{x}=1\}.\]
It is easy to see that $\norm{\cdot}^\ast$ satisfies the norm axioms.
For any $x\in V$ there exists a functional $\fhi\in V^\ast$ with
\[ \abs{\fhi(x)} = \norm{\fhi}^\ast\norm{x},\]
which we abbreviate as $x\perp_{\norm{\cdot}}\fhi$, or $x\perp\fhi$ if the norm is clear from the context.
(This is the Hahn-Banach theorem, which is geometrically obvious in dimension $2$: any convex body has a supporting line at each boundary point.)
A unit $\fhi$ for which $x\perp\fhi$ is usually called a \emph{norming functional} of $x$.

\subsection{Bilinear forms}
An \emph{identification} of $V$ and $V^\ast$ is an isomorphism $T:V\to V^\ast$.
Although there is no ``natural'' identification of a finite-dimensional vector space and its dual, we will now attempt to find a natural as possible identification.
Any identification of $V$ and $V^\ast$ corresponds to a bilinear form on $V$, defined by $[x,y]:=T(y)(x)$.
It is easily seen that, since $T$ is an isomorphism, $[\cdot,\cdot]$ will be non-degenerate: if $[a,y]=0$ for all $y\in V$, then $a=o$ (equivalently, if $[x,b]=0$ for all $x\in V$, then $x=o$).
It is also easy to see that conversely, any non-degenerate bilinear form $[\cdot,\cdot]$ defines an isomorphism $T:V\to V^\ast$ by $T(x)=[x,\cdot]$.
\begin{lemma}
A bilinear form on a vector space satisfies
\begin{equation}\label{orth}
 [x,y]=0 \iff [y,x]=0 \quad\text{for all $x,y\in V$}
\end{equation}
if and only if either
\begin{itemize}
\item $[x,x]=0$ for all $x\in V$ \textup{(}in which case we call it a \emph{symplectic form}\textup{)}, or
\item $[x,y]=[y,x]$ for all $x,y\in V$ \textup{(}in which case we call it an \emph{orthogonal form} or \emph{inner product}\textup{)}.
\end{itemize}
\end{lemma}
Thus, if we make assumption \eqref{orth}, there are two types of bilinear forms to choose from.
For a proof, see Artin \cite[Chapter~III.2]{MR90h:51003}.
The next lemma is an easy exercise.
\begin{lemma}
Up to a non-zero multiple, there is only one non-degenerate symplectic form on a two-dimensional vector space.
If coordinates are chosen, then a symplectic form is a non-zero multiple of the determinant of the $2\times2$ matrix which has the two vectors as columns.
\end{lemma}
Thus a choice of a symplectic form on $V$ is the same as a choice of area unit plus orientation.
By contrast, there are three degrees of freedom in choosing an orthogonal form on a two-dimensional $V$.
It can be shown \cite[Chapter~III.7]{MR90h:51003} that each orthogonal form in the plane corresponds uniquely to the curve $[x,x]=1$, which is either an ellipse (and then we have Euclidean geometry), or a hyperbola (and then we have two-dimensional spacetime, with one space and one time dimension).

It is therefore much more natural to use a symplectic form in two-dimensional affine geometry (when we do not want to impose a Euclidean structure).
This is indeed also what we need to define the antinorm.
\section{The antinorm}
From now on we assume that the plane $V$ has norm $\norm{\cdot}$ and symplectic bilinear form $[\cdot,\cdot]$.
The \emph{antinorm} on $V$ is now defined to be the dual norm on $V^\ast$ identified with $V$ via $[\cdot,\cdot]$, i.e.\
\[ \antinorm{x}:=\norm{Tx}^\ast=\sup\{[y,x]:\norm{y}=1\},\]
and we obtain
\begin{equation}\label{one}
\antinorm{x} = \sup\{[x,y]:\norm{y}=1\}.
\end{equation}
It also follows that
\begin{equation}\label{onea}
\abs{[x,y]}\leq\norm{x}\antinorm{y}.
\end{equation}
We denote the unit ball of the antinorm by
\[I := \{v\in V: \antinorm{v}\leq 1\},\]
and the unit \emph{anticircle} by
\[\bd I := \{v\in V: \antinorm{v}= 1\}.\]
The unit anticircle is also called the \emph{isoperimetrix}, and was introduced by Busemann \cite{MR9:372h}.
In the literature it is usually defined by choosing a Euclidean structure in the plane.
Then the isoperimetrix is the polar body of the unit ball, rotated by an angle of $90^\circ$.
The above definition is simpler.

\begin{proposition}
The antinorm of the antinorm \textup{(}still with respect to the same bilinear form\textup{)} is the original norm.
\end{proposition}
\begin{proof}
Let $\norm{x}_\mathrm{a,a}$ denote the antinorm of the antinorm.
Then by taking the supremum of \eqref{onea} over all $y$ with $\antinorm{y}=1$ and using \eqref{one} we obtain $\norm{x}_\mathrm{a,a}\leq\norm{x}$.

Secondly we let $\fhi$ be a norming functional of $x$, i.e.\ $\norm{\fhi}=1$ and $\fhi(x)=\norm{x}$.
Then $\fhi=Ty$ for some $y\in V$ with $\antinorm{y}=1$.
Thus $[x,y]=\fhi(x)=\norm{x}$, and \eqref{one} gives $\norm{x}_\mathrm{a,a}\geq[x,y]=\norm{x}$.
\end{proof}
It follows that
\begin{equation}\label{two}
x\perp_{\norm{\cdot}}Ty\iff\abs{[x,y]}=\norm{x}\antinorm{y}\iff y\perp_{\antinorm{\cdot}}Tx.
\end{equation}

\section{Normality}
\begin{figure}[h]
\begin{center}
\begin{overpic}[scale=0.5]{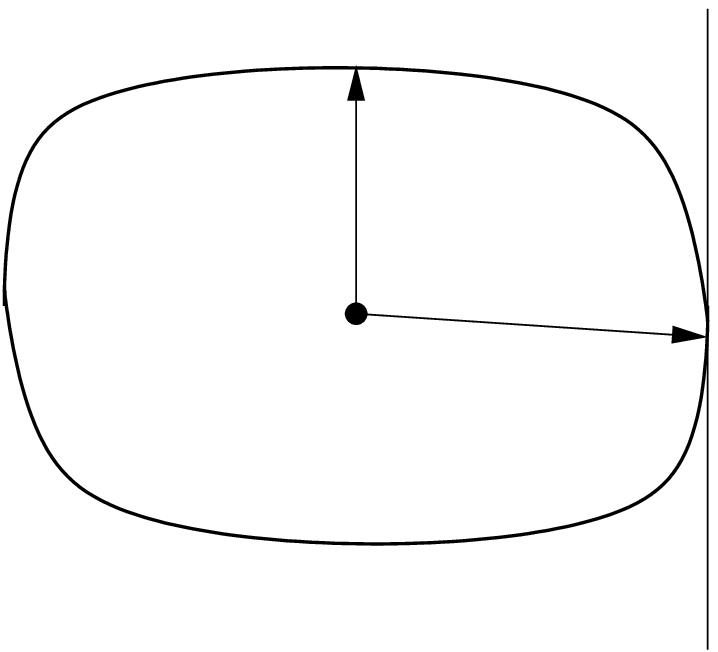}
\put(90,36){$x$}
\put(55,75){$y$}
\end{overpic}
\end{center}
\caption{$x\normal y$}\label{fig1}
\end{figure}

The relation of \emph{normality} was introduced by Carath\'eodory \cite{Blaschke, Radon}.
It is usually referred to as Birkhoff orthogonality, due to its rediscovery in \cite{Birkhoff}.
A non-zero $x\in V$ is \emph{normal} to a non-zero $y\in V$, denoted $x\normal y$, if
\[ \norm{x}\leq\norm{x+\lambda y}\quad\text{for all $\lambda\in\R$}.\]
Geometrically this means that the line through $\frac{1}{\norm{x}}x$ parallel to $y$ is a supporting line of the unit ball at $\frac{1}{\norm{x}}x$ (Figure~\ref{fig1}).

It follows immediately from \eqref{one} that 
\begin{equation}\label{three}
x\perp_{\norm{\cdot}} Ty \iff x\normal y,
\end{equation}
and that
\begin{equation}\label{four}
x\normal y\iff x\text{ maximizes $\abs{[\cdot,y]}$ on $V$.}
\end{equation}

\begin{theorem}[Busemann \cite{MR9:372h}]\label{thm1}
The antinorm reverses normality and is the unique such norm up to a choice of unit area.
\end{theorem}
\begin{proof}
The first part is immediate from \eqref{two} and \eqref{three}.
Uniqueness follows from the fact that if two norms have the same normality relation, then the norms are proportional.
This follows immediately from Lemma~\ref{tech1}.
For a detailed proof see Sch\"opf \cite{MR98m:46018}.
\end{proof}
We denote normality with respect to the antinorm by $x\antinormal y$.
This relation, opposite to that of normality, is also called \emph{transversality} \cite[Section~4.6]{MR97f:52001}.


\section{Minkowski content}
Given a Minkowski plane in which $[\cdot,\cdot]$ has been fixed, the \emph{Minkowski content} of a segment $ab$ is defined to be
\[ \mnorm{ab} := \lim_{\epsi\to 0^+}\frac{\area{ab+\epsi B}}{2\epsi}.\]
This notion was introduced by Minkowski \cite{Minkowski2}.
Here \[ab+\epsi B := \{x+y : x\in ab, y\in \epsi B\}\] is the \emph{Minkowski sum} of the sets $ab$ and $\epsi B$, i.e., the union of all balls of radius $\epsi$ centred at a point of $ab$.
It is easily seen that this limit exists: Choose a unit $y$ such that $y\normal a-b$.
Then
\begin{eqnarray*}
\area{ab+\epsi B}&=&\abs{[a-b,y]}+\area{\epsi B}\\
&=& 2\epsi\antinorm{a-b}+\epsi^2\area{B} \quad\text{(by \eqref{two} and \eqref{three}),}
\end{eqnarray*}
and
\[ \frac{\area{ab+\epsi B}}{2\epsi} = \antinorm{a-b}+\frac\epsi2\area{B}.\]
Thus $\mnorm{ab}=\antinorm{a-b}$, which gives
\begin{theorem}
The Minkowski content of a segment is its length in the antinorm.
\end{theorem}

\section{Radon curves}
In general it is not true that $\normal$ is a symmetric relation.
Radon \cite{Radon} constructed norms for which $\normal$ is symmetric and proved that his construction gives all such norms.
The unit circles of these norms are called \emph{Radon curves}.
Later Birkhoff \cite{Birkhoff} and Day \cite{MR9:192c} also gave constructions.
These constructions are all formulated in terms of polarity and a $90^\circ$ rotation with respect to some auxiliary Euclidean structure.
We present the construction using only the bilinear form $[\cdot,\cdot]$.

\begin{figure}[h]
\begin{center}
\bigskip\bigskip
\begin{overpic}[scale=0.5]{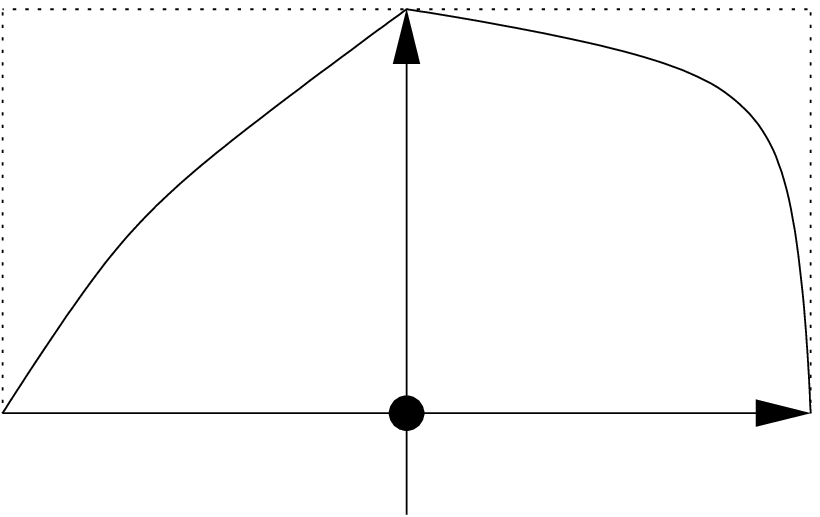}
\put(93.5,6){$a$}
\put(53,54){$b$}
\put(55,18){$Q_1$}
\put(35,18){$Q_2$}
\put(35,5){$Q_3$}
\put(55,5){$Q_4$}
\put(86,55){$C_1$}
\put(20,50){$C_2$}
\end{overpic}
\end{center}
\caption{Constructing a Radon curve}\label{fig2}
\end{figure}

Choose any linearly independent $a, b\in V$ with $[a,b]=1$.
If we define a coordinate system with $a$ and $b$ as the standard unit vectors, then $[a,b]$ is just the determinant, and we may speak of the four quadrants $Q_1, Q_2, Q_3, Q_4$.
See Figure~\ref{fig2}.
In $Q_1$ choose any curve $C_1$ joining $a$ and $b$ such that $C_1$ is contained in the parallelogram with vertices $o,a,b,a+b$, and is on the boundary of $B_1:=\conv(C_1\cup\{o\})$.
This piece of a unit ball $B_1$ already defines a norm in the first quadrant:
\[ \norm{x} := \inf\{\lambda:\lambda^{-1}x\in B_1\} \quad\text{for all $x\in Q_1$}.\]
We want to extend $C_1$ to a curve $C_2$ in $Q_2$ in such a way that normality is symmetric.
(Once the unit circle has been determined in $Q_2$, it is fixed in $Q_3$ and $Q_4$ by central symmetry.)
Note that the direction vectors of supporting lines of $C_1$ ($C_2$), correctly chosen, all lie in $Q_2$ ($Q_1$, respectively).
Thus what we want is that $x\normal y\iff y\normal x$ for all $x\in C_1,y\in C_2$.
Note that we can already use the partially defined norm in $Q_1$ to define the antinorm in $Q_2$:
\[ \antinorm{x} = \sup\{\abs{[x,y]}:y\in Q_1,\norm{y}=1\}\quad\text{for all $x\in Q_2$.}\]
It follows from what we already know about antinorms that $x\normal y\iff y\normal_{\antinorm{\cdot}}x$ for all $x\in C_1, y\in Q_2$.
It follows that the norm defined by $C_2$ must have the same normality relation as the antinorm, so it must be a multiple of the antinorm.
Thus $C_2$ must be a multiple of the unit anticircle in $Q_2$.
Due to convexity, $C_2$ must join $b$ and $-a$, so we find that $C_2$ is exactly the unit anticircle in $Q_2$.
The construction is finished.

This construction gives a large class of Minkowski planes for which normality is symmetric.
Also for any unit vectors $x,y$ with $x\normal y$ we have $\abs{[x,y]}=1$.
By scaling we then obtain further norms such that $x\normal y\implies\abs{[x,y]}=\lambda$ for some fixed $\lambda>0$ independent of $x$ and $y$.
All such norms are called \emph{Radon norms}, and their unit circles \emph{Radon curves}.
If $\abs{[x,y]}=1$ for all unit $x,y$ for which $x\normal y$, then we say that the Radon norm and Radon curve are \emph{normalized}.

Conversely, if we are given a norm $\norm{\cdot}$ on $V$ for which normality is symmetric, we choose any unit $a,b$ with $a\normal b$ (and $b\normal a$).
We scale $[\cdot,\cdot]$ by some $\lambda>0$ so that we may assume $[a,b]=1$.
Define quadrants as before using $a$ and $b$.
We take $C_1$ to be the unit circle in $Q_1$, and do the construction as before to obtain a Radon norm $\norm{\cdot}'$ which then has the same normality relation as $\norm{\cdot}$.
Since their unit circles already coincide in the first quadrant, we obtain $\norm{\cdot}=\norm{\cdot}'$, which means that the given norm is a Radon norm.
This gives
\begin{theorem}[Radon \cite{Radon}]\label{radonthm}
A unit circle in a Minkowski plane is a Radon curve if and only if normality is symmetric.
\end{theorem}
The above discussion also shows that in a Radon norm with $[a,b]=1$ for some $a,b$ with $a\normal b$, the antinorm coincides with the norm.
So we also have
\begin{corollary}[Busemann \cite{MR9:372h}]\label{antinorm}
A norm is Radon if and only if it equals a multiple of its antinorm.
\end{corollary}
The following equivalent statement was already given by Radon.
\begin{corollary}[Radon \cite{Radon}]
The Minkowski content of a segment coincides with its length in the norm if and only if the plane is Radon.
\end{corollary}
Thus for Radon norms there is a natural choice of area unit, namely that for which the antinorm equals the norm.
This happens exactly when the Radon curve is normalized.
Thus for Radon norms we have a true natural identification of the plane and its dual as for inner product spaces.

As examples of Radon curves we mention ellipses (which is also the historical origin of the term conjugate diameters) and regular $(4n+2)$-gons.
Also the mixed $\ell_p$-$\ell_q$ norm on $\R^2$ is a Radon norm for any $1<p,q<\infty$ satisfying the conjugacy relation $p^{-1}+q^{-1}=1$:
\[ \norm{(\alpha,\beta)}_{p,q} := \begin{cases} (\abs{\alpha}^p+\abs{\beta}^p)^{1/p} & \text{if $\alpha\beta\geq 0$},\\ (\abs{\alpha}^q+\abs{\beta}^q)^{1/q} & \text{if $\alpha\beta\leq 0$}.\end{cases}.\]

Gruber proved the following stability versions of Theorem~\ref{radonthm} and its corollaries.
Intuitively they say that if normality is almost symmetric, or if the norm is almost a multiple of its antinorm, then the norm is almost a Radon norm.
We say that normality is \emph{symmetric up to} $\epsi>0$ if for any unit $x,y$ such that $x\normal y$ there exists a unit $z$ with $y\normal z$ and $\norm{x-z}\leq\epsi$.
\begin{theorem}[Gruber \cite{Gruber3}]
There exists a constant $\gamma>0$ with the following property.
Let $(V,\norm{\cdot})$ be a Minkowski plane and $\epsi>0$ be given satisfying one of the following two properties:
\begin{itemize}
\item normality is symmetric up to $\epsi>0$,
\item for some $\lambda>0$ we have that the norm and the antinorm scaled by $\lambda$ differ on the unit sphere by at most $\epsi$:
\[  \bigabs{\norm{x}-\lambda\antinorm{x}}\leq\epsi\;\text{for all}\; x \;\text{with}\; \norm{x}=1.\]
\end{itemize}
Then there exists a Radon norm $\norm{\cdot}'$ such that the \emph{Banach-Mazur distance} between $\norm{\cdot}$ and $\norm{\cdot}'$ is at most $1+\gamma\epsi$:
\[ \norm{x}\leq\norm{x'}\leq(1+\gamma\epsi)\norm{x} \text{ for all $x\in V$.}\]
\end{theorem}

James \cite{MR9:42c} studied the normality relation in the context of Banach spaces.
One of his results for Radon norms in fact gives a characterization of Radon curves (Corollary~\ref{3prime}).
In Theorem~\ref{3} we give the corresponding result valid for all norms.

\begin{theorem}\label{3}
For any $x,y\neq o$, if $x\normal \lambda x+y$ and $\mu y+x\normal y$, then $0\leq \lambda\mu\leq 2$.
\end{theorem}
\begin{proof}
%
%
Without loss of generality $\lambda\mu\neq 1$.
From the definition of $\normal$ it follows in particular that
\[ \norm{x+\tfrac{\mu}{1-\lambda\mu}(\lambda x+y)}\geq\norm{x} \text{ and } \norm{\mu y+x-\mu y} \geq\norm{\mu u+x}.\]
Simplifying we obtain
\[ \norm{x+\mu y}\geq\abs{1-\lambda\mu}\,\norm{x} \text{ and } \norm{x}\geq\norm{\mu y+x}.\]
It follows that $\abs{1-\lambda\mu}\leq1$.
\end{proof}
We need the following lemma that will also be used later.
\begin{lemma}[James \cite{MR9:42c}]\label{james}
For any linearly independent $x,y\neq o$ there exists an $\alpha\in\R$ such that $x\normal\alpha x+y$.
Furthermore, for any such $\alpha$ we have $\abs{\alpha}\leq\norm{y}/\norm{x}$.
\end{lemma}
\begin{proof}
Choose $z\neq o$ such that $x\normal z$.
Then $z=\lambda x+\mu y$ for some $\lambda,\mu\in\R$.
If $\mu=0$ then $\alpha\neq 0$ and $x\normal x$, a contradiction.
Therefore $x\normal \mu^{-1}\lambda x+y$, and we take $\alpha=\mu^{-1}\lambda $.

By the definition of normality we then get for all $\lambda\in\R$ that $\norm{x}\leq\norm{x+\lambda(\alpha x+y)}=\norm{-(1/\alpha) y}$, if we take $\lambda=-1/\alpha$.
The inequality follows.
\end{proof}
The forward implication of the following corollary of Theorem~\ref{3} is from James \cite{MR9:42c}.
The other direction was stated without proof in \cite[Proposition~38]{MR2002h:46015}.
\begin{corollary}[James \cite{MR9:42c}, M-S-W \cite{MR2002h:46015}]\label{3prime}
A norm is Radon if and only if the following holds:
For any $x,y\neq o$, if $x\normal \lambda x+y$ and $y\normal \mu y+x$, then $\lambda\mu\geq 0$.
\end{corollary}
\begin{proof}
The $\Rightarrow$ direction is immediate from Theorem~\ref{3}.

$\Leftarrow$: Suppose that for some unit $x, y$ we have $x\normal y$ but $y\notnormal x$.
By Lemma~\ref{james} we have $y\normal x+\lambda_0 y$ for some $\lambda_0\neq 0$.
Without loss of generality $\lambda_0>0$ (replacing $y$ by $-y$, if necessary).
Since $y\notnormal x$ it follows by continuity that for all sufficiently small $\alpha\in(0,\lambda_0^{-1})$ there exist $\beta\in(0,\lambda_0)$ such that $y+\alpha x\normal x+\beta y$.
Thus $0<\alpha\beta<1$, and letting $y'=y+\alpha x$ we obtain $y'\normal x+\beta(1-\alpha\beta)^{-1}y'$.
On the other hand, $x\normal y'-\alpha x$.
Since $-\alpha<0$, it follows from the hypothesis that $\beta(1-\alpha\beta)^{-1}\leq 0$.
Thus $1-\alpha\beta<0$, a contradiction.
\end{proof}

\section{Triangle area}
The $\frac12 \mathrm{base}\times\mathrm{height}$ formula for the area of a triangle in Euclidean geometry generalizes as follows.
Given any triangle $\triangle a_1a_2a_3$, we let $\beta_i$ be the length in the norm of the side opposite $a_i$, and $\eta_i$ the the shortest distance from $a_i$ to the line through the side opposite $a_i$ measured in the norm, and call it the \emph{height}.
Similarly, we let $\eta_{\mathrm{a}i}$ be the the shortest distance from $a_i$ to the line through the side opposite $a_i$ measured in the antinorm, and call it the \emph{anti-height}.
\begin{proposition}\label{trianglearea}
The area of $\triangle a_1a_2a_3$ is $\frac12 \beta_i\eta_{\mathrm{a}i}$.
\end{proposition}
\begin{proof}
Let $v=a_3-a_2$, choose $p$ on the line $a_2a_3$ such that $\antinorm{a_1-p}=\eta_{\mathrm{a}1}$.
(Thus $p$ is a ``foot of the perpendicular from $a_1$''.)
Let $u=a_1-p$.
Then $u\antinormal v$, hence $v\normal u$.
It follows that \[\area{\triangle a_1a_2a_3}=\frac12\bigabs{[u,v]}=\frac12\norm{v}\,\antinorm{u}\]
by \eqref{two} and \eqref{three}.
\end{proof}
Averkov \cite[Theorem~5.1]{Averkov} gives a generalization of this formula.
In the next characterization of Radon curves, the $\Rightarrow$ direction was observed by Busemann \cite{MR17:779a}, and also by \golab\ according to Tam\'assy \cite{MR23:A4052}.
The $\Leftarrow$ direction was proved by Tam\'assy \cite{MR23:A4052} in the smooth case, but in the more general context of starshaped unit circles that are not necessarily centrally symmetric.
\begin{corollary}[\golab-Busemann \cite{MR23:A4052, MR17:779a} $\Rightarrow$, Tam\'assy \cite{MR23:A4052} $\Leftarrow$]\label{gbt}
A norm is Radon if and only if for all triangles $\triangle a_1a_2a_3$ the value $\frac12 \beta_i\eta_i$ is independent of $i$.
\end{corollary}
\begin{proof}
For the $\Rightarrow$ direction note that if the norm is Radon, $\eta_{\mathrm{a}i}=\lambda\eta_i$ for some fixed $\lambda$, by Corollary~\ref{antinorm}.

$\Leftarrow$: Let $x$ and $y$ be unit vectors with $x\normal y$.
The height of $\triangle oxy$ from $x$ is then $\norm{x}=1$.
Let $\eta$ be the height of $\triangle oxy$ from $y$.
By hypothesis, $\frac12\norm{y}\norm{x}=\frac12\norm{x}\eta$.
See Figure~\ref{fig3}.
\begin{figure}[h]
\begin{center}
\bigskip
\begin{overpic}[scale=0.5]{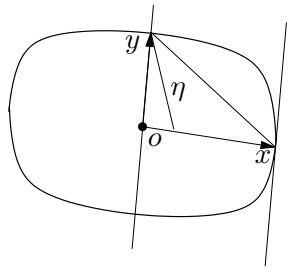}
\put(88,37){$x$}
\put(50,43){$o$}
\put(42,78){$y$}
\put(58,61){$\eta$}
\end{overpic}
\end{center}
\caption{Proof of Corollary~\ref{gbt}}\label{fig3}
\end{figure}
Thus $\eta=\norm{y}$, which means that $\norm{y-o}$ is also the shortest distance from $y$ to the line $ox$.
Thus the line $ox$ supports the unit ball with centre $y$.
Equivalently, the line through $y$ parallel to $ox$ supports $B$, i.e., $y\normal x$.
By Theorem~\ref{radonthm} the norm is Radon.
\end{proof}
Averkov \cite[Theorem~5.2]{Averkov} gives a related characterization of Radon norms.
For further results on area, see Section~\ref{polygonarea}.

\section{Angles and bisectors}
\subsection{Angular bisectors}\label{bisectors}
In Euclidean geometry an angular bisector has the following two well-known characterizations in terms of distance:
\begin{itemize}
\item Any point on it has the same distance to the two sides of the angle,
\item In $\triangle abc$, if $ad$ is a bisector of $\myangle a$, with $d$ on $bc$, then $bd/dc=ba/ac$.
\end{itemize}
We may use any of these two properties to extend the notion of angular bisector to Minkowski planes.

\begin{figure}[h]
\begin{center}
\includegraphics[scale=0.5]{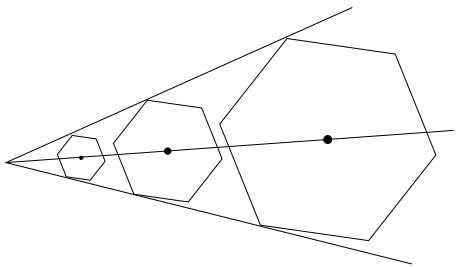}
\end{center}
\caption{Glogovskii angular bisector}\label{fig4}
\end{figure}

Glogovskii \cite{MR45:2597} uses the first property.
It is easily seen that the points equidistant (in the norm) to the two sides of an angle lie on a line.
This line is called the \emph{Glogovskii angular bisector} of the angle.
See Figure~\ref{fig4}.

Busemann \cite{MR51:8961} uses the second property.
It is again easily seen that for any angle $\myangle a$, the set of points $d$ for which there exists a line through $d$ intersecting the two sides of $\myangle a$ in $b$ and $c$, say, such that $\norm{b-d}/\norm{d-c}=\norm{b-a}/\norm{a-c}$, also lies on a line (and then all lines through $d$ intersecting the sides of $\myangle a$ satisfy this property).
This line is called the \emph{Busemann angular bisector} of the angle.
A simple way of constructing the Busemann angular bisector is to choose two points $b$ and $c$ on the two sides of $\myangle a$ with $\norm{b-a}=\norm{c-a}$, and then to let $d$ be the midpoint of $b$ and $c$.
The following proposition follows immediately from the definition:
\begin{proposition}
The Glogovskii bisectors of the three angles of a triangle are concurrent, and the point of intersection is the centre of the \textup{(}unique\textup{)} inscribed Minkowski circle of the triangle.
\end{proposition}
It can be proved by Ceva's theorem that the Busemann angular bisectors of the angles of a triangle are also concurrent.
A more conceptual explanation of this fact follows from the following theorem.
\begin{theorem}[D\"uvelmeyer \cite{Nico-cmb}]\label{angletheorem}
The Busemann angular bisector in the norm coincides with the Glogovskii angular bisector in the antinorm \textup{(}and vice versa\textup{)}.
\end{theorem}
\begin{proof}
Consider any angle $\myangle bac$.
Assume without loss of generality that $\norm{b-a}=\norm{c-a}$.
Let $d$ be the midpoint of $b$ and $c$.
Then, as observed above, $od$ is the Busemann bisector of $\myangle bac$.
Triangles $\triangle abd$ and $\triangle acd$ clearly have the same area, and their bases $ab$ and $ac$ have the same length in the norm.
By Proposition~\ref{trianglearea} they have the same anti-height, i.e., the distance in the antinorm from $d$ to the line $ab$ equals that to the line $ac$.
Thus $od$ is the Glogovskii bisector in the antinorm.
\end{proof}
\begin{corollary}\label{busemannbi}
The Busemann bisectors of the three angles of a triangle are concurrent, and the point of intersection is the centre of the \textup{(}unique\textup{)} inscribed anticircle of the triangle.
\end{corollary}
\begin{corollary}[D\"uvelmeyer \cite{Nico-jg}]
A norm is Radon if and only if Busemann and Glogovskii angular bisectors coincide.
\end{corollary}
\begin{proof}
The $\Rightarrow$ direction is clear from Theorem~\ref{angletheorem}.

Conversely, consider the unit ball $B$ with centre $o$.
Let $\ell_0$ be some fixed supporting line of $B$.
Let $\lambda I+v$ be the anticircle such that $\ell_0$ supports $\lambda I+v$ as well.
Thus clearly $-\ell_0$ also supports $\lambda I+v$.

Let $\ell$ be any other supporting line of $B$ not parallel to $\ell_0$.
Then $\ell$ and $\ell_0$ determines an angle with vertex $a$, say.
Thus the Glogovskii bisector of this angle is $ao$.
By the hypothesis this is also the Busemann bisector, so it is also the Glogovskii bisector in the antinorm, by Theorem~\ref{angletheorem}.
It follows that $\ell$ also supports $\lambda I+v$.

By Lemma~\ref{tech2} then $B=\lambda I+v$, i.e., the norm is a multiple of the antinorm, which gives that the norm is Radon by Theorem~\ref{antinorm}.
\end{proof}

An \emph{angular measure} is a Borel measure $\mu$ on the unit circle $\bd B$ of a Minkowski plane such that
\begin{itemize}
\item $\mu(\bd B)=2\pi$, 
\item for any Borel $S \subset \bd B$ we have $\mu(S)=\mu(-S)$, and 
\item for each $p\in\bd B$ we have $\mu(\{p\})=0$.
\end{itemize}
The measure of an angle $\myangle a$ is then defined to be the measure of the arc of the unit circle determined by the angle translated to the origin.
Brass \cite{MR97c:52036} used a special such angular measure to analyze packings of unit circles in Minkowski planes.
For any angular measure satisfying the following additional property one may define a corresponding angular bisector:
\begin{itemize}
\item All non-degenerate arcs of $\bd B$ have positive measure.
\end{itemize}
D\"uvelmeyer \cite{Nico-cmb} proved that if the Busemann or Glogovskii bisector can be defined using such an angular measure, then the norm must be Euclidean.

Note that there is a unique angular measure $\mu_l$ such that $\mu_l(\myangle o)$ is proportional to the length of the arc of the unit circle determined by $\myangle o$.
There is also a unique angular measure $\mu_a$ such that $\mu_a(\myangle o)$ is proportional to the area of the sector of the unit circle determined by $\myangle o$.
\begin{theorem}
In any Minkowski plane, $\mu_a(\myangle o)$ is proportional to the length in the antinorm of the arc of the unit circle determined by $\myangle o$.
\end{theorem}
\begin{proof}
This is immediate from Proposition~\ref{trianglearea} if the unit ball is a polygon.
In general we obtain the result by approximating the unit ball by centrally symmetric polygons.
We omit the details.
\end{proof}
The dual of this (with norm and antinorm interchanged) is called a ``Kepler law'' by Wallen \cite{Wallen}:
\begin{quote}
\emph{``If we travel with constant speed along an isoperimetrix, equal areas are swept out \textup{(}from the center\textup{)} in equal times.''}
\end{quote}
It follows that in a Radon plane $\mu_a=\mu_l$.
However, there are other norms for which they coincide as well, the most obvious example being when the unit circle is a square.
Helfenstein \cite{Helfenstein} asked as a problem to determine the norms for which these two angular measures coincide.
Unfortunately his solution \cite{Helfenstein2} claimed that this happens exactly when the norm is Radon.
D\"uvelmeyer \cite{Nico-cmb} correctly characterized these norms.
They are exactly those whose unit circle is an \emph{equiframed curve}, i.e., a centrally symmetric convex curve for which each point is touched by a circumscribed parallelogram of minimum area.
(See \cite{MS} for more on equiframed curves.)

\subsection{Perpendicular bisectors}
Finally, we mention two results on the generalization of Euclidean perpendicular bisectors.
The \emph{bisector} of two points $p$ and $q$ in a Minkowski plane is the set of points equidistant from $p$ and $q$:
\[B(p,q) := \{x\in V: \norm{x-p}=\norm{x-q}\}.\]
Bisectors in Minkowski spaces have been studied mostly by computational geometers; see \cite{MSII} for a survey.
It is known that bisectors are lines if and only if the plane is Euclidean (Mann \cite{Mann}).
If the Minkowski plane is \emph{strictly convex} (i.e.\ $\norm{x+y}<\norm{x}+\norm{y}$ for linearly independent $x,y$), then bisectors are always unbounded curves.
On the other hand, if the Minkowski plane is not strictly convex, it is not difficult to find a bisector containing interior points.
Thus a Minkowski plane is strictly convex if and only if all bisectors are curves (see \cite{MR2002h:46015, MSII}).
The same considerations also give that a Minkowski plane is strictly convex if and only if each bisector is contained in some strip bounded by two parallel lines (see \cite{MSII} for references).
The relation to the antinorm is as follows.
\begin{theorem}\label{bisector}
In a strictly convex plane, the bisector of $x$ and $y$ is contained in the interior of a unique strip, which has the property that its bounding lines are tangent to the anticircle with diameter $xy$.
\end{theorem}
\begin{proof}
Since the plane is strictly convex, its anticircles are smooth, i.e., all their supporting lines are tangent lines.
Choose $v\neq o$ such that $v\normal y-x$.
Since then $y-x\antinormal v$, we have that the lines tangent at $x$ and $y$ to the anticircle with diameter $xy$ must be parallel to $v$.
See Figure~\ref{fig5}.

\begin{figure}[h]
\begin{center}
\bigskip
\begin{overpic}[scale=0.5]{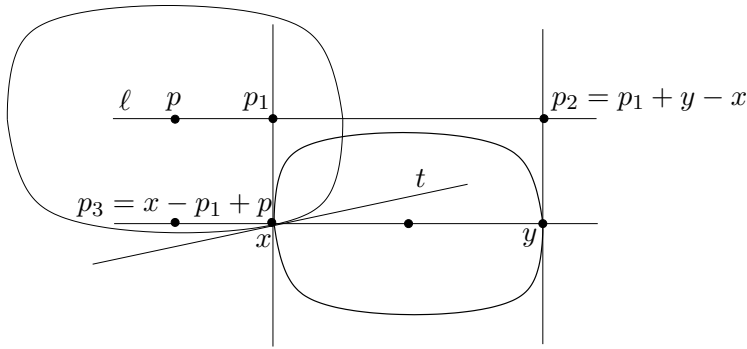}
\put(19,40){$\ell$}
\put(27,41){$p$}
\put(40,41){$p_1$}
\put(92,41){$p_2=p_1+y-x$}
\put(12,23){$p_3=x-p_1+p$}
\put(42,17){$x$}
\put(69,27){$t$}
\put(87,18){$y$}
\end{overpic}
\end{center}
\caption{Proof of Theorem~\ref{bisector}}\label{fig5}
\end{figure}
Consider an arbitrary line $\ell$ parallel to $xy$ intersecting the tangent lines in $p_1$ and $p_2=p_1+y-x$, say.
It is sufficient to show that points on $\ell$ not between $p_1$ and $p_2$ are not equidistant to $x$ and $y$.
Let $p$ be a point on $\ell$ such that $p_1$ is between $p$ and $p_2$.
Since $v\normal y-x$, we have that the shortest distance from $x$ to $\ell$ is attained by $\norm{x-p_1}$.
Because of strict convexity, no other point on $\ell$ attains the shortest distance, hence $\norm{p-x}>\norm{p_1-x}$.
Let $p_3=x-p_1+p$.
Thus $\norm{p-p_3}<\norm{p-x}$.
Consider the ball with centre $p$ and radius $px$.
Thus $p_3$ is in its interior.
Let $t$ be a supporting line of this ball at $x$.
Then $y$ and $p_3$ are on opposite sides of $t$, and in particular, $y$ is not in the ball, i.e., $\norm{p-y}>\norm{p-x}$.

Similarly, if $p\in\ell$ is such that $p_2$ is between $p$ and $p_1$, then $\norm{p-y}<\norm{p-x}$.
\end{proof}
\begin{corollary}\label{radonbisector}
A Minkowski plane is strictly convex and Radon if and only if for all $x,y$ the bisector of $x,y$ is contained in the strip whose bounding lines support the circle with diameter $x,y$.
\end{corollary}
\begin{proof}
The $\Rightarrow$ direction is immediate from the previous theorem.

For the converse choose any $a,b$ with $a\normal b$.
Then the bisector $B(o,a)$ is contained in two strips, one parallel to $v$ where $a\antinormal v$ (by the previous theorem), and one parallel to $b$ (by hypothesis).
Since $B(x,y)$ is easily seen to be unbounded, we must have that $v$ and $b$ are parallel.
It follows that $a\antinormal b$, i.e., $b\normal a$, and the result follows from Theorem~\ref{radonthm}.
\end{proof}

\section{Anti-equilateral triangles}
A triangle in a Minkowski plane is \emph{anti-equilateral} if it is equilateral in the antinorm.
Anti-equilateral triangles appear in the Fermat-Torricelli problem for triangles, and in the problem of reduced convex bodies in the plane.

The following is an immediate corollary of Proposition~\ref{trianglearea}.
\begin{corollary}[Viviani's theorem]
The sum of the distances from any point inside an anti-equilateral triangle to its sides is a constant, where the distance of a point inside the triangle is considered to be positive, and outside the triangle to be negative.
\end{corollary}
\begin{proposition}\label{incentre}
A triangle is anti-equilateral if and only if its incentre coincides with its centroid.
\end{proposition}
\begin{proof}
Let $p$ be the incentre and $s$ the centroid of $\triangle abc$.
Asume first that $\triangle abc$ is anti-equilateral.
Then, since $\triangle sab$, $\triangle sbc$, and $\triangle sac$ have the same area, and $ab$, $bc$, and $ac$ are the same length in the antinorm, the distance (in the norm) from $s$ to the three sides of $\triangle abc$ must be equal (by Proposition~\ref{trianglearea}), i.e., $s$ must be an incentre.
Since the incentre is unique, $s=p$.

The converse follows along similar lines.
\end{proof}

A point $p$ is a Fermat-Torricelli point of the triangle $\triangle abc$ in a Minkowski plane if $p$ minimizes the function $x\mapsto\norm{x-a}+\norm{x-b}+\norm{x-c}$, i.e., if $p$ is a point such that the sum of its distances to the vertices of the triangle is a minimum.

\begin{theorem}
If $p\neq a,b,c$, then $p$ is a Fermat-Torricelli point of $\triangle abc$ if and only if $p$ satisfies the following property:
\begin{quote}
Let $a',b',c'$ be the points of intersection of the rays $\Ray{pa}$, $\Ray{pb}$, $\Ray{pc}$ with a circle centred at $p$.
Then there exist supporting lines to the circle at $a',b',c'$ forming an anti-equi\-la\-te\-ral triangle.
\end{quote}
\end{theorem}
\begin{proof}
The above statement, with the words ``an anti-equilateral triangle'' replaced by ``a triangle with centroid $p$'', was proved in \cite{MSW}, and in the special case of smooth norms, in \cite{CG}.
The theorem now follows from Proposition~\ref{incentre}.
\end{proof}
The $\Leftarrow$ direction of the above theorem can also be proved by adapting Viviani's proof of the characterization of the Fermat-Torricelli point of a triangle in the Euclidean plane (see D\"orrie \cite[Problem~91]{Dorrie}).
See \cite{MSW} for more on the Fermat-Torricelli problem in Minkowski spaces.

The \emph{minimum width} of a convex body $C$ is the minimum distance (in the norm) between two parallel supporting lines of $C$, where the minimum is taken over all pairs of parallel supporting lines of $C$.
A convex body is \emph{reduced} if it does not properly contain a convex body of the same minimum width.
See \cite{MR92b:52016} and \cite{LM} for more on reduced convex bodies in Minkowski spaces.
Averkov \cite[Theorem~5.3]{Averkov} found the following characterization of reduced triangles in a Minkowski plane.
\begin{theorem}
A triangle is reduced in a Minkowski plane if and only if it is anti-equilateral.
\end{theorem}

\section{Area and perimeter of polygons and convex curves}\label{polygonarea}
Another theorem from Euclidean geometry is the following.
\begin{quote}
For any convex polygon $P$ circumscribed about a circle of radius $\rho$ we have
\[ \rho\perimeter{P}=2\area{P},\]
where $\perimeter{P}$ is the perimeter of $P$.
\end{quote}
It has a very simple proof, and as can be guessed from Proposition~\ref{trianglearea}, the corresponding statement for Minkowski planes must involve the anticircle.
\begin{theorem}\label{circump}
For any convex polygon $P$ containing an anticircle of radius $\rho$ we have
\[ \rho\perimeter{P}\leq2\area{P}\]
with equality if and only if $P$ is circumscribed about the anticircle.
In particular, all polygons of the same area circumscribed about an anticircle have the same perimeter.
\end{theorem}
\begin{proof}
Immediate from Proposition~\ref{trianglearea}, if we subdivide $P$ into triangles by joining the vertices of $P$ to the centre of the anticircle.
\end{proof}
Using Lemma~\ref{tech3} to approximate a convex body by circumscribed polygons, we obtain in the limit the following result, which in the Euclidean case is mentioned in \cite[p.~82, eq.~10]{Zbl.8:77}.
The second statement is mentioned for example by Ghandehari \cite{Ghandehari}.
\begin{corollary}\label{circumpcor}
For any convex body $C$ containing an anticircle of radius $\rho$ we have
\[ \rho\perimeter{C}\leq2\area{C}.\]
Equality holds in particular for $C$ an anticircle \textup{(}but also for other convex bodies\textup{)}:
The perimeter of an anticircle $\rho I$ satisfies
\[ \rho\perimeter{\rho I} = 2\area{\rho I}.\]
\end{corollary}
We now present a new characterization of Radon curves based on the above corollary.
Note that its proof uses the results on angular bisectors discussed above in Section~\ref{bisectors}.
\begin{corollary}\label{newcor}
The following statements are equivalent for a Minkowski plane:
\begin{enumerate}
\item\label{i}The Minkowski plane is Radon,
\item\label{ii} for some fixed $\gamma>0$ we have that for any convex polygon $P$ circumscribed about a circle of radius $\rho$, $\perimeter{P}=\gamma\rho\area{P}$,
\item\label{iii} for some fixed $\gamma>0$ and some fixed $n\geq 3$, we have that for any convex $n$-gon $P$ circumscribed about a circle of radius $\rho$, $\perimeter{P}=\gamma\rho\area{P}$,
\item\label{iv} for some fixed $\gamma>0$ we have that for any triangle $P$ circumscribed about a circle of radius $\rho$, $\perimeter{P}=\gamma\rho\area{P}$.
\end{enumerate}
\end{corollary}
\begin{proof}
\ref{i}$\implies$\ref{ii} follows from the previous theorem and Corollary~\ref{antinorm}.

\ref{ii}$\implies$\ref{iii} is trivial.

\ref{iii}$\implies$\ref{iv} follows since we may approximate circumscribed triangles by circumscribed $n$-gons; see Figure~\ref{fig6}.
\begin{figure}[h]
\begin{center}
\includegraphics[scale=0.5]{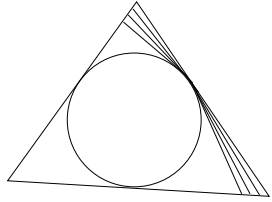}
\end{center}
\caption{Approximating a triangle by an $n$-gon}\label{fig6}
\end{figure}

\ref{iv}$\implies$\ref{i}: Let $x$ and $y$ be unit vectors with $x\normal y$ and $y\normal x$.
Denote the four quadrants defined by $x,y$ by $Q_1,\dots,Q_4$ as before.
See Figure~\ref{fig7}.
\begin{figure}[h]
\begin{center}
\begin{overpic}[scale=0.4]{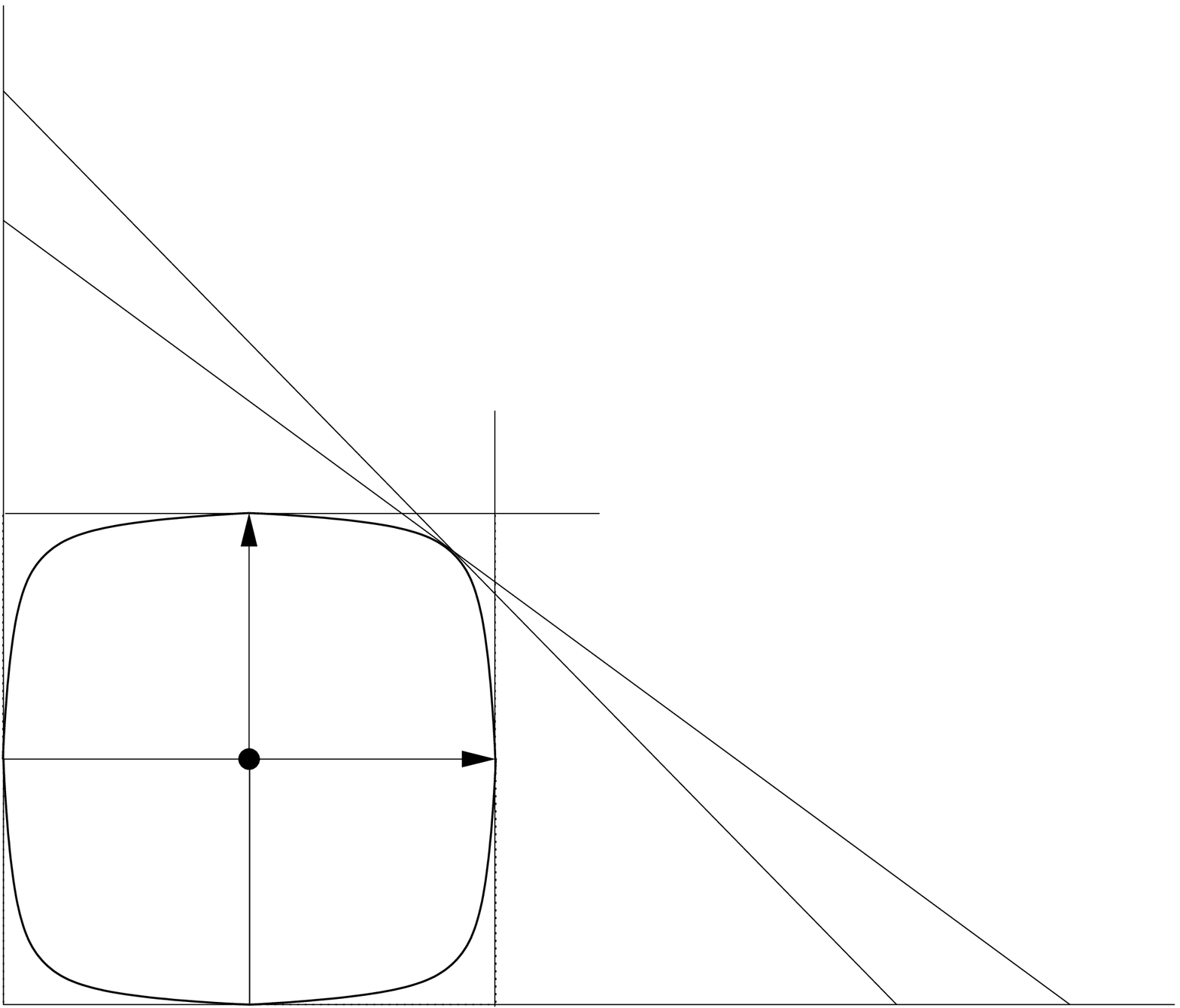}
\put(-8,48){$-\ell_y$}
\put(43,48){$\ell_y$}
\put(23,39){$y$}
\put(52,41){$\ell_x$}
\put(33,35){$B$}
\put(23,25){$Q_1$}
\put(14,25){$Q_2$}
\put(14,16.5){$Q_3$}
\put(23,16.5){$Q_4$}
\put(70,17){$\ell_1$}
\put(60,10){$\ell_0$}
\put(53,2){$-\ell_x$}
\end{overpic}
\end{center}
\caption{Proof of Corollary~\ref{newcor}}\label{fig7}
\end{figure}
Note that $x$ and $y$ have the same antinorm, say $\rho:=\antinorm{x}=\antinorm{y}$.
We show that $\rho I=B$ by showing that all supporting lines of $B$ also support $\rho I$.

We already have this for the lines $\ell_x$ parallel to $x$ passing through $y$, and $\ell_y$ parallel to $y$ passing through $x$.
Choose a fixed line $\ell_0\neq\ell_x,\ell_y$ supporting $B$ in the first quadrant $Q_1$.
Let $\ell_1\neq\ell_x,\ell_y$ be any other line supporting $B$ in $Q_1$.
Then $-\ell_x$, $-\ell_y$, $\ell_i$ determines a triangle $\Delta_i$, for $i=0,1$, both circumscribed to $B$.
By hypothesis,
\begin{equation}\label{x}
\frac{\area{\Delta_0}}{\perimeter{\Delta_0}} =\frac{\area{\Delta_1}}{\perimeter{\Delta_1}}.
\end{equation}
Let $I_i$ be the anticircle inscribed in $\Delta_i$.
The centres of both anticircles are on the Busemann bisector of the angle determined by $-\ell_x$ and $-\ell_y$ (by Corollary~\ref{busemannbi}).
Also, by Theorem~\ref{circump} and \eqref{x} they have the same radius (in the antinorm) $\lambda$, say.
Thus they coincide, say $I_0=I_1=:I'$.
It follows that $I'$ shares all supporting lines of $B$ in $Q_1$.
By taking limits it follows that $I'$ also has $\ell_x$ and $\ell_y$ as supporting lines.
Since now $\ell_x$ and $-\ell_x$ are both supporting lines of $I'$, it follows that $I'=\rho I$.
Thus $\rho I$ has the same supporting lines as $B$ in $Q_1$.
A similar argument gives that $\rho I$ and $B$ share the same supporting lines in $Q_2$.
By central symmetry, $\rho I$ and $B$ also share the same supporting lines in $Q_3$ and $Q_4$.
By Lemma~\ref{tech2}, $\rho I=B$, and by Corollary~\ref{antinorm} the norm is Radon.
\end{proof}

Finally we mention the following very interesting duality.
\begin{theorem}[Sch\"affer \cite{MR47:5732}, Thompson \cite{MR52:4138}]
The perimeter in the norm of the unit anticircle equals the perimeter in the antinorm of the unit circle.
\end{theorem}

\section{Isoperimetric inequalities and the Zenodorus problem}
The \emph{isoperimetric problem} is the problem to find, among all closed rectifiable curves of a fixed length in the norm, those of largest area.
It is easily seen that we may restrict our attention to convex curves, since taking the convex hull of a curve does not increase the length and does not decrease the area.
Similarly, the \emph{Zenodorus problem} is the problem of finding, among all $n$-gons of a fixed perimeter, those of largest area.
If we define the \emph{isoperimetric ratio} of a convex body $C$ (with respect to a given norm) to be
\[ \iota(C) = \frac{\perimeter{C}^2}{\area{C}},\]
then the isoperimetric problem is to find all $C$ of smallest $\iota(C)$, and the Zenodorus problem is to find, for each $n$, all $n$-gons $P$ of smallest $\iota(P)$, where the minimum is taken over all $n$-gons.

L.~Fejes T\'oth \cite[\S4]{FT} approaches both problems in the Euclidean plane in a discrete way by first proving an isoperimetric inequality for polygons (see \eqref{*} below).
One may then approximate an arbitrary convex body by circumscribed polygons and take the limit of the inequality, to obtain an inequality of Bonnesen \cite[p.~82, eq.~11]{Zbl.8:77} (see \eqref{+} below), which solves the isoperimetric problem in the Euclidean plane without needing to assume the existence of a solution.
Fejes T\'oth notes that this approach, using inner parallel polygons, comes from Sz.~Nagy \cite{SzNagy}.
See also Niven \cite[Ch.~12]{Niven} for a careful exposition.

The same may be done in Minkowski planes, as shown by Chakerian \cite{Chakerian}.
However, Chakerian only states a corollary of \eqref{*} which generalizes Lhuilier's inequality in Euclidean geometry (see \eqref{**} below and \cite[\S4]{FT}).
By examining his proof one sees that (after taking limits) he really generalizes the inequality of Bonnesen referred to above.
Thus, implicit in his paper is a complete solution to the isoperimetric problem, analogous to that of Sz.~Nagy mentioned above.
Secondly, as shown below, using his approach one also obtains a solution to the Zenodorus problem for Minkowski planes.
This seems not to have been mentioned before in the literature.

\begin{theorem}[implicit in Chakerian \cite{Chakerian}]\label{chakerian}
Let $P$ be a convex polygon and let $\rho I$ be the largest anticircle contained in $P$.
Let $Q$ be the polygon circumscribed about $\rho I$ with sides parallel to those of $P$.
Then
\begin{equation}\label{*}
\area{P}+\area{Q}\leq \rho\perimeter{P}.
\end{equation}
\end{theorem}
See \cite{Chakerian} for the proof.
\begin{corollary}\label{chakeriancor}
\begin{equation}\label{--}
\iota(P)-\iota(Q)\geq \frac{1}{\area{P}}(\perimeter{P}-\perimeter{Q})^2.
\end{equation}
\end{corollary}
\begin{proof}
Inequality \eqref{*} is algebraically equivalent to
\begin{equation}\label{-}
\frac{\perimeter{P}^2}{\area{P}} - 4\frac{\area{Q}}{\rho^2}\geq\frac{1}{\area{P}}\left(\perimeter{P}-2\frac{\area{Q}}{\rho}\right)^2,
\end{equation}
Then use Theorem~\ref{circump} to rewrite \eqref{-} as \eqref{--}.
\end{proof}
This implies the Minkowski equivalent of Lhuilier's inequality.
\begin{corollary}[Chakerian \cite{Chakerian}]
\begin{equation}\label{**}
\perimeter{P}^2\rho^2\geq 4\area{P}\area{Q},
\end{equation}
with equality if and only if $P=Q$, i.e., if $P$ is circumscribed about an anticircle.
\end{corollary}
We now present our solution to the Zenodorus problem.
We define a \emph{Zenodorus $n$-gon} to be a convex $n$-gon $P$ circumscribed about an anticircle $I'$, such that it has the smallest area among all $n$-gons circumscribed about $I'$.
By a simple compactness argument, there exists at least one Zenodorus $n$-gon for each $n\geq 3$.
(It is possible for certain norms that there is, up to scaling, only one Zenodorus $n$-gon for a fixed $n$.)
\begin{corollary}\label{zenodorus}
Among all $n$-gons of a fixed perimeter the ones of largest area are the Zenodorus $n$-gons.
\end{corollary}
\begin{proof}
As mentioned above, the solutions to the Zenodorus problem must be convex.
Thus let $P$ be a convex $n$-gon of smallest $\iota(P)$.
Let $\rho I$ and $Q$ be as in Theorem~\ref{chakerian}.
By choice of $P$, $\iota(Q)\geq\iota(P)$.
Combining this with \eqref{--} (Corollary~\ref{chakeriancor}), we obtain $\perimeter{P}=\perimeter{Q}$ and $\iota(P)=\iota(Q)$.
Thus $\area{P}=\area{Q}$ and $P=Q$.
Thus all $P$ of smallest $\iota(P)$ are circumscribed about an anticircle.
By scaling we may assume this to be the unit anticircle $I$.
Then by Theorem~\ref{circump}, $\perimeter{P}=2\area{P}$, hence $\iota(P)=4\area{P}$.
Thus, in order to minimize $\iota(P)$, we only have to minimize $\area{P}$ among all $P$ circumscribed about the unit anticircle.
\end{proof}
Thus the Zenodorus problem for Minkowski planes is reduced to a problem in convex geometry that does not refer to distances.
In general it is not possible to say much about $n$-gons of smallest area circumscribed about a convex body.
However, the following two general statements can be made.
\begin{proposition}[Day \cite{MR9:246h}, Dowker \cite{MR5:153m}]
Let $P_n$ be a convex $n$-gon of smallest area circumscribed about a convex body $C$.
\begin{description}
\item The sides of $P_n$ touch $C$ at their midpoints.
\item The areas of $P_n$ form a convex sequence: \[\area{P_{n-1}}+\area{P_{n+1}}\geq2\area{P_n}.\]
\end{description}
\end{proposition}
Using Lemma~\ref{tech3} to approximate a convex body $C$ by polygons and taking the limit in Theorem~\ref{chakerian} we obtain the Minkowski generalization of the inequality of Bonnesen mentioned above:
\begin{corollary}\label{thm14}
Let $C$ be a planar convex body and let $\rho I$ be the largest anticircle contained in $C$.
Then
\begin{equation}\label{+}
\area{C}+\area{\rho I}\leq\rho\perimeter{C}.
\end{equation}
\end{corollary}
Wallen \cite{Wallen} has a similar inequality.
In the same way as in the proof of Corollary~\ref{chakeriancor} (using Corollary~\ref{circumpcor} instead of Theorem~\ref{circump}), we obtain
\[ \frac{\perimeter{C}^2}{\area{C}} - \frac{\perimeter{I}^2}{\area{I}} \geq 0\]
with equality if and only if $\rho\perimeter{C}=2\area{\rho I}$, which together with \eqref{+} is equivalent to $\area{C}=\area{\rho I}$.
Thus we have obtained
\begin{theorem}[Busemann \cite{MR9:372h}]\label{isoperimetricthm}
The only figures solving the isoperimetric problem in a Minkowski plane are the anticircles.
\end{theorem}
Note that we did not need the existence of a convex body attaining $\min\iota(C)$; this is a corollary of the proof.
By Corollary~\ref{antinorm} we then have
\begin{corollary}[Busemann \cite{MR9:372h}]
A norm is Radon if and only if circles solve the isoperimetric problem.
\end{corollary}
Corollary~\ref{thm14} has the following analogue for smallest anticircles containing $C$.
\begin{theorem}\label{thm15}
Let $C$ be a planar convex body and let $\sigma I$ be the smallest anticircle containing $C$.
Then
\begin{equation}\label{++}
\area{C}+\area{\sigma I}\leq\sigma\perimeter{C}.
\end{equation}
\end{theorem}
Inequalities \eqref{+} and \eqref{++} together are algebraically equivalent to
\begin{equation}\label{+++}
\area{C}-\alpha\perimeter{C}+\alpha^2\area{I}\leq 0 \text{ for all $\alpha\in[\rho,\sigma]$.}
\end{equation}
This inequality is a special case of an inequality of Blaschke on mixed area (in our case $\perimeter{C}$ is the mixed area of $C$ and $I$) \cite{Blaschke}.
Also, \eqref{+++} implies the following extension of Bonnesen's isoperimetric deficit inequality (as observed by Petty \cite{MR18:760e}):
\[ \perimeter{C}^2 - 4\area{I}\area{C}\geq\area{I}^2(\sigma-\rho)^2. \]
Unfortunately, we do not know an elementary proof of \eqref{++} involving polygons as in the proof of \eqref{+}.
A proof of Blaschke's inequality (which includes \eqref{++}) may be found in \cite{Flanders}, and sharpenings may be found in \cite{MR94k:52003}.

\section{Non-expansive mappings}
In the Euclidean plane it is elementary that for any point $x$ on a circle and any line $\ell$ through the centre of the circle, the foot of the perpendicular from $x$ on $\ell$ is inside the circle.
This statement is false for general norms.
To change it into a true statement it is however sufficient to consider nearest points in the antinorm.
\begin{lemma}\label{nearestpoint}
Let $x$ and $y$ be unit vectors in a Minkowski plane.
Then any point on the line $oy$ nearest to $x$ in the antinorm is in the unit ball $B$.
\end{lemma}
\begin{proof}
Let $m$ be a point on $oy$ nearest to $x$ in the antinorm.
Clearly $m-x\antinormal y$.
Thus $y\normal m-x$ by Theorem~\ref{thm1}.
If $m=\lambda y$, then $y\normal \lambda y-x$, and by Lemma~\ref{james}, $\norm{m}=\abs{\lambda}\leq\norm{x}/\norm{y}=1$.
\end{proof}
\begin{corollary}[De Figueiredo and Karlovitz \cite{MR35:2130}, Amir {\cite[\S 18]{MR88m:46001}}]
The following statements are equivalent for a Minkowski plane:
\begin{enumerate}
\item The plane is Radon.
\item Let $x$ and $y$ be any two unit vectors $x,y$.
Then any point on the line $oy$ nearest to $x$ in the norm satisfies $p\in B$.
\item For any unit vectors $x, y$ there exists a point $p$ on the line $oy$ nearest to $x$ in the norm, satisfying $p\in B$.
\end{enumerate}
\end{corollary}
A proof is contained in \cite[Proposition~40]{MR2002h:46015}.

The \emph{radial projection} of a Minkowski plane onto its unit ball is defined by
\[ p(x) := \begin{cases}
x & \text{if $\norm{x}\leq 1$,}\\
\frac{1}{\norm{x}}x & \text{if $\norm{x}>1$.}
\end{cases} \]
A mapping $f:V\to V$ is \emph{non-expansive} if $\norm{f(v)-f(w)}\leq\norm{v-w}$ for all $v,w\in V$.
It is easily seen that in the Euclidean plane the radial projection is non-expansive.
By now the following theorem comes without surprise.
\begin{theorem}[Karlovitz \cite{MR46:7869}]\label{nonexpansive}
The radial projection is non-expansive in the antinorm.
\end{theorem}
\begin{proof}
Let $v,w\in V$.
We have to show that $\antinorm{f(v)-f(w)}\leq\antinorm{v-w}$.
If $\norm{v},\norm{w}\leq 1$, then this is trivial.
If $\norm{w}\geq\norm{v}\geq 1$, then $\antinorm{w-v}\geq\antinorm{\frac{1}{\norm{v}}w-\frac{1}{\norm{v}}v}$, so in this case it is sufficient to consider the case $\norm{w}\geq\norm{v}=1$.

Thus the only case that we have to consider is $\norm{w}\geq1\geq\norm{v}$.
Let $p$ be a point on $ow$ nearest to $x$ in the antinorm.
By Lemma~\ref{nearestpoint} $p$ is in the ball with centre $o$ and radius $\norm{x}$, i.e., $\norm{p}\leq\norm{x}\leq 1$.
Thus, if we write $w=\lambda u$ and $p=\mu u$, with $\lambda=\norm{w}$, then we have $\mu\leq 1\leq \lambda$, and $\antinorm{\mu u-v}\leq\antinorm{\lambda u-v}$.
Applying the triangle inequality we obtain
\begin{align*}
\antinorm{\frac{1}{\norm{w}}w-v} & = \antinorm{\frac{\lambda -1}{\lambda-\mu}(\mu u-v) + \frac{1-\mu}{\lambda-\mu}(\lambda u-v)}\\
&\leq \frac{\lambda -1}{\lambda-\mu}\antinorm{\mu u-v} + \frac{1-\mu}{\lambda-\mu}\antinorm{\lambda u-v}\\
&\leq \frac{\lambda -1}{\lambda-\mu}\antinorm{\lambda u-v} + \frac{1-\mu}{\lambda-\mu}\antinorm{\lambda u-v}\\
&\leq \antinorm{\lambda u - v}=\antinorm{w-v}.\qedhere
\end{align*}
\end{proof}

\begin{corollary}[De Figueiredo and Karlovitz \cite{MR35:2130}]
A norm is Radon if and only if the radial projection is non-expansive in the norm.
\end{corollary}
A proof is in \cite[Proposition~40]{MR2002h:46015}.

One may define the radial projection of any set $S$ that is \emph{starshaped with respect to the origin}, i.e., such that the intersection of $S$ and any line through the origin is a segment.
By interchanging norm and antinorm in Theorem~\ref{nonexpansive} we obtain that the radial projection onto an antiball with centre $o$ is non-expansive in the norm.
The antiballs are unique with respect to this property:
\begin{theorem}[Gruber \cite{MR57:1279}]
Let $S$ be starshaped with repsect to the origin in a Minkowski plane.
Then the radial projection onto $S$ is non-expansive if and only if $S$ is an antiball with centre $o$.
\end{theorem}

\section{Generalized convexity and abstract approximation theory}
The following three types of generalized convexity are defined only in terms of the metric, and therefore provide a way of studying convexity in arbitrary metric spaces.
We consider them only in Minkowski planes.
\subsubsection*{$d$-convexity}
This is called \emph{Minkowski convexity} by Petty \cite{MR18:760e}.
We first need to define the notion of a $d$-segment, introduced by Menger \cite{Me}.
The \emph{$d$-segment} of $a$ and $b$ is the set of all points \emph{metrically between} $a$ and $b$:
\[ [a,b]_d = \{x\in V: \norm{a-b}=\norm{a-x}+\norm{x-b}\}.\]
A set $S$ is then called \emph{$d$-convex} if it contains $[a,b]_d$ for all $a,b\in S$.

Since $[a,b]_d$ contains the ordinary segment $ab$, $d$-convex sets are convex.
In strictly convex normed spaces, $d$-segments are always ordinary segments, and then all convex sets are $d$-convex.
However, in general it is possible for convex sets to be not $d$-convex.
We may think of $d$-convexity as a type of ``superconvexity''.

\subsubsection*{$\norm{\cdot}$-convexity}
This was introduced by Menger \cite{Me}, who called it \emph{metric convexity}.
We first need to define the notion of a metric segment.
Any rectifiable curve joining $a$ and $b$ has a length in the norm $\norm{\cdot}$, which is always at least $\norm{a-b}$.
A curve joining $a$ and $b$ of length exactly $\norm{a-b}$ is called a \emph{metric segment} between $a$ and $b$.
A set $S$ is then called \emph{$\norm{\cdot}$-convex} if for any $a,b\in S$ some metric segment between $a$ and $b$ is contained in $S$.

Since ordinary segments are metric segments, convex sets are $\norm{\cdot}$-convex.
In  strictly convex normed spaces, metric segments are always ordinary segments, and then all $\norm{\cdot}$-convex sets are convex.
However, in general it is possible for $\norm{\cdot}$-convex sets to be non-convex.
We may think of $\norm{\cdot}$-convexity as a type of ``subconvexity''.

\subsubsection*{$B$-convexity}
This notion of convexity was introduced by Lassak \cite{MR57:1271}, and studied by him and others in many further papers; see the survey \cite{MSII} for further references.
A set $S$ is $B$-convex if for any finite subset $A$ of $S$, the intersection of all balls containing $A$ is contained in $S$.
Since balls are convex, it follows that $B$-convex sets are convex.
In smooth normed spaces the intersection of all balls containing a finite set $A$ is the ordinary convex hull of $A$, and it follows that then convex sets are $B$-convex.
However, in general it is possible for convex sets to be not $B$-convex.
We may also think of $B$-convexity as a type of ``superconvexity''.

The following result connects $B$-convexity and $d$-convexity in Minkowski planes.
\begin{theorem}[Lassak \cite{MR85a:52005}]
A set is $d$-convex in the norm if and only if it is $B$-convex in the antinorm.
\end{theorem}
Since balls are obviously $B$-convex, the above result gives that antiballs are $d$-convex.
This also follows from the isoperimetric property of anticircles (Theorem~\ref{isoperimetricthm}), as noted by Petty \cite{MR18:760e}.
\begin{corollary}
In a Radon plane $d$-convex and $B$-convex sets coincide.
\end{corollary}
We do not obtain a characterization of Radon planes, since e.g.\ for all strictly convex and smooth planes $d$-convexity and $B$-convexity coincide with the usual convexity.

We now consider certain notions from abstract approximation theory.
Given a closed set $S$, the distance from any point $x$ to $S$ is defined as
\[ d(x,S) = \inf\{d(x,s):s\in S\}.\]
The \emph{metric projection onto $S$} is the set-valued mapping
\[ p_S(x) = \{s\in S: \norm{x-s}=d(x,S)\},\quad x\in V.\]
Because $S$ is closed, $p_S(x)$ will always be non-empty.
The set $S$ is a \emph{Chebyshev set} if $p_S(x)$ is a singleton for all $x\in V$.
In this case we consider $p_S$ to be an ordinary function from $V$ to $S$.
The set $S$ is also called a \emph{C2 set} if $p_S(x)$ is a contractible set.
Thus closed convex sets are always C2 sets.
In a Minkowski plane, $p_S(x)$ is always a subset of a circle.
Thus $S$ is a C2 set if $p_S(x)$ is an arc.
The connection to $\norm{\cdot}$-convexity is as follows.
\begin{theorem}[Gruber \cite{MR86k:41032}]
A set is C2 if and only if it is closed and $\antinorm{\cdot}$-convex.
\end{theorem}
Gruber also characterizes Chebyshev sets using the related notion of semistrict $\antinorm{\cdot}$-convexity; see \cite{MR86k:41032} for the definition.
Also see Gruber \cite{MR57:1279} for related results.

We remark that Hetzelt \cite{MR86h:41037} considers certain related notions involving ``co-ap\-proxi\-ma\-tion'' and uses the antinorm.

The theorems on non-expansiveness of radial projections may be generalized as follows:
\begin{theorem}[Karlovitz \cite{MR46:7869}]
Let $S$ be a convex Chebyshev set in a Minkowski plane.
Then its metric projection is non-expansive in the antinorm.
\end{theorem}
The result of Karlovitz is in fact more general: he extends the definition of metric projection to any closed convex set in a Minkowski plane, not necessarily one that is Chebyshev.
This also extends the radial projection of the unit ball defined above.

\begin{corollary}[Phelps \cite{MR19:432a}]
A norm is Radon if and only if for any convex Chebyshev set its metric projection is non-expansive in the norm.
\end{corollary}
Phelps only considers the strictly convex case, but it can be seen that this restriction is not necessary.

\section{Higher dimensions}
Normality between vectors can still be defined in higher dimensions, although it now makes more sense to define the normality relation between a vector and a hyperplane.
It is known \cite{Blaschke} that if normality between vectors is symmetric and the dimension is at least $3$, then the Minkowski space must be Euclidean.
It also follows that if all two-dimensional subspaces of a Minkowski space of dimension at least $3$ are Radon, then the norm must again be Euclidean.
Many of the above characterizations of Radon curves then lead to an analogous characterization of Euclidean spaces of dimension at least $3$.
See Amir \cite{MR88m:46001} for an extensive list of characterizations of Euclidean spaces.

The antinorm can, as before, be defined as the Minkowski content of a segment.
This gives a norm, but it does not seem to be known how the unit ball of the Minkowski content is related to the original unit ball.
According to Gruber \cite{Gruber2, Gruber3} it is conjectured that the Minkowski content and the norm are proportional only for Euclidean spaces.

The isoperimetric problem also has higher-dimensional analogues.
Here there are different solutions depending on how $(n-1)$-dimensional measure (``area'') is defined; see \cite{MR97f:52001}.

In conclusion we note that since the norm/antinorm duality depends on an antisymmetric bilinear form, which is unique up to a scalar multiple only for a two-dimensional space, the phenomena discussed in this paper are in a certain sense essentially two-dimensional.
One may speculate however that for certain normed space of even dimension symplectic forms may also bring out certain Euclidean features.


\begin{thebibliography}{88}

\bibitem{MR88m:46001}
D.~Amir, \emph{Characterizations of inner product spaces}, Birkh\"auser
  Verlag, Basel, 1986.

\bibitem{MR90h:51003}
E.~Artin, \emph{Geometric Algebra}, John Wiley \& Sons Inc., New York, NY, 1988. Reprint of the 1957 original.



\bibitem{Averkov}
G. Averkov, \emph{On the geometry of simplices in Minkowski spaces}, Studies of the University of \v{Z}ilina, Mathematical Series \textbf{14} (2001), 1--13.


\bibitem{Birkhoff}
G.~Birkhoff, \emph{Orthogonality in linear metric spaces}, Duke Math.\ J.
  \textbf{1} (1935), 169--172.

\bibitem{Blaschke}
W.~Blaschke, \emph{{R\"a}umliche {V}ariationsprobleme mit symmetrischer
  {T}ransversalit{\"a}tsbedingung}, Ber.\ Verh.\ S{\"a}chs.\ Ges.\ Wiss.\ Leipzig.\
  Math.-Phys.\ Kl.\ \textbf{68} (1916), 50--55.

\bibitem{Zbl.8:77}
T. Bonnesen, W. Fenchel: \emph{Theorie der
konvexen K\"orper}. Ergebn.\ Math.\ und ihrer Grenzgebiete, Bd.~3,
Springer, Berlin 1934.

\bibitem{MR97c:52036}
P.~Brass, \emph{Erd{\H o}s distance problems in normed spaces}, Comput. Geom.
  \textbf{6} (1996), 195--214.

\bibitem{MR9:372h}
H.~Busemann, \emph{The isoperimetric problem in the {M}inkowski plane}, Amer.\ J.
  Math.\ \textbf{69} (1947), 863--871.

\bibitem{MR17:779a}
H.~Busemann, \emph{The {G}eometry of {G}eodesics}, Academic Press Inc., New York,
  N.Y, 1955.

\bibitem{MR51:8961}
H. Busemann: \emph{Planes with analogues to Euclidean angular bisectors},
Math.\ Scand.\ \textbf{36} (1975), 5--11.

\bibitem{Chakerian}
G.~D.~Chakerian, \emph{The isoperimetric problem in the Minkowski plane}, Amer.\ Math.\ Monthly \textbf{67} (1960), 1002--1004.

\bibitem{CG}
G.~D. Chakerian and M.~A. Ghandehari, \emph{The {Fermat} problem in {Minkowski}
  spaces}, Geom. Dedicata \textbf{17} (1985), 227--238.

\bibitem{Coxeter}
H.~S.~M.~Coxeter, Introduction to Euclidean Geometry, 2nd ed., Wiley, 1989.

\bibitem{MR9:246h}
M.~M. Day, \emph{Polygons circumscribed about closed convex curves}, Trans.\
  Amer.\ Math.\ Soc.\ \textbf{62} (1947), 315--319.

\bibitem{MR9:192c}
M.~M.~Day, \emph{Some characterizations of inner-product spaces}, Trans.\ Amer.\
  Math.\ Soc.\ \textbf{62} (1947), 320--337.

\bibitem{Dorrie}
H.~D\"orrie, \emph{100 Great Problems of Elementary Mathematics: Their History and Solution}, translated by David Antin, Dover, New York, 1965.
Original: \emph{Triumph der Mathematik}, Physica-Verlag, W\"urzburg, 1965.

\bibitem{MR5:153m}
C.~H. Dowker, \emph{On minimum circumscribed polygons}, Bull. Amer. Math. Soc.
  \textbf{50} (1944), 120--122.

\bibitem{Nico-jg}
N.~D\"uvelmeyer, \emph{A new characterization of Radon curves via angular bisectors}, to appear in Journal of Geometry.

\bibitem{Nico-cmb}
N.~D\"uvelmeyer, \emph{Angle measures and bisectors in Minkowski planes}, to appear in Canadian Mathematical Bulletin.

\bibitem{FT}
L.~Fejes T\'oth, \emph{Lagerungen in der Ebene auf der Kugel und im Raum}, 2nd ed., Springer, Berlin, 1972.

\bibitem{MR35:2130}
D.~G.~de~Figueiredo and L.~A. Karlovitz, \emph{On the radial projection in
  normed spaces}, Bull.\ Amer.\ Math.\ Soc.\ \textbf{73} (1967), 364--368.

\bibitem{Flanders}
H.~Flanders, \emph{A Proof of Minkowski's Inequality for Convex Curves}, Amer.\ Math.\ Monthly \textbf{75} (1968), 581--593.

\bibitem{Ghandehari}
M.~A.~Ghandehari, \emph{Steinhardt's inequality in the Minkowski plane}, Bull.\ Austral.\ Math.\ Soc.\ \textbf{45} (1992), 261--266.

\bibitem{MR45:2597}
V.~V. Glogovs'ki{\u\i}, \emph{Bisectors on the {M}inkowski plane with norm
  $(x\sp{p}+y\sp{p})\sp{1/p}$ ({R}ussian)}, V\=\i snik L'v\=\i v. Pol\=\i tehn.
  \=Inst. (1970), 192--198, 218.

\bibitem{MR57:1279}
P.~M.~Gruber: \emph{Fixpunktmengen von Kontraktionen in endlich-dimensionalen 
nor\-mier\-ten R\"aumen}, Geom.\ Dedicata \textbf{4} (1975), 179--198.

\bibitem{MR86k:41032}
P.~M.~Gruber, \emph{Planar {C}hebyshev sets}, Mathematical Structures -- Computational
  Mathematics -- Mathematical Modelling, 2, Bulgar.\ Acad. Sci., Sofia, 1984, 
  pp.~184--191.

\bibitem{Gruber2}
P.~M.~Gruber, \emph{Radons Beitr\"age zur Konvexit\"at}, In: Johann Radon, Collected Works, Vol.~1, Verlag der \"Osterreichischen Akademie der Wissenschaften, Vienna, 1987.

\bibitem{Gruber3}
P.~M.~Gruber, \emph{Stability of Blaschke's characterization of ellipsoids and Radon norms}, Discrete Comp.\ Geom.\ \textbf{17} (1997), 411--427.


\bibitem{Helfenstein}
H.~Helfenstein, \emph{Problem 13}, Canad.\ Math.\ Bull.\ \textbf{2} (1959), 43.

\bibitem{Helfenstein2}
H.~Helfenstein, \emph{Solution to Problem 13}, Canad.\ Math.\ Bull.\ \textbf{4} (1961), 77--78.

\bibitem{MR86h:41037}
L.~Hetzelt: \emph{On suns and cosuns in finite-dimensional normed real vector spaces}, Acta Math.\ Hungar.\ \textbf{45} (1985), 53--68.

\bibitem{MR9:42c}
R.~C.~James, \emph{Orthogonality and linear functionals in normed linear spaces},
  Trans.\ Amer.\ Math.\ Soc.\ \textbf{61} (1947), 265--292.

\bibitem{MR46:7869}
L.~A.~Karlovitz, \emph{The construction and application of contractive
  retractions in $2$-dimensional normed linear spaces}, Indiana Univ.\ Math.\ J.
  \textbf{22} (1972/73), 473--481.

\bibitem{MR57:1271}
M. Lassak: \emph{On metric $B$-convexity for which diameters of any set
and its hull are equal}, Bull.\ Acad.\ Polon.\ Sci.\ S\'er.\ Sci.\ Math.\ Astronom.\ 
Phys.\ \textbf{25} (1977), 969--975.

\bibitem{MR85a:52005}
M.~Lassak, \emph{Some connections between $B$-convexity and
$d$-convexity}, Demonstratio Math.\ \textbf{15} (1982),
261--270.

\bibitem{MR92b:52016}
M.~Lassak, \emph{Reduced convex bodies in the plane}, Israel J. Math \textbf{70} (1990), 365--379.

\bibitem{LM} M.~Lassak and H.~Martini, \emph{Reduced bodies in Minkowski space}, to appear in Acta Math.\ Hungarica.

\bibitem{Mann}
H.~Mann: \emph{Untersuchungen \"uber Wabenzellen bei allgemeiner
Minkowskischer Metrik}, Monatsh.\ Math.\ Phys.\ \textbf{42}
(1935), 417--424.

\bibitem{MS}
H.~Martini, K.~J.~Swanepoel, \emph{Equiframed curves --- a generalization of Radon curves},
Monatsh.\ Math.\ \textbf{141} (2004), 301--314.

\bibitem{MSII}
H.~Martini, K.~J.~Swanepoel, \emph{The geometry of Minkowski spaces --
a survey. Part II}, Expo.\ Math.\ \textbf{22} (2004), 93--144.

\bibitem{MR2002h:46015}
H.~Martini, K.~J.~Swanepoel, G.~Weiss, \emph{The geometry of Minkowski spaces --
a survey. Part I}, Expo.\ Math.\ \textbf{19} (2001), 97--142.
(Errata: Expo.\ Math.\ \textbf{19} (2001), p.~364.)

\bibitem{MSW}
H. Martini, K.~J.~Swanepoel, G. Weiss: \emph{The Fermat-Torricelli
problem in normed planes and spaces},  J. Optim.\ Theory Appl.\ \textbf{115} (2002), 
283--314.

\bibitem{Me}
K. Menger: \emph{Untersuchungen \"uber allgemeine Metrik},
Math.\ Ann.\ \textbf{100} (1928), 75--163.

\bibitem{Minkowski}
H.~Minkowski, \emph{Geometrie der {Z}ahlen}, B. G. Teubner, Leipzig und Berlin,
  1896, 1910.

\bibitem{Minkowski2}
H.~Minkowski, \emph{\"Uber die Begriffe L\"ange, Oberfl\"ache und Volumen}, Jahresber.\ Deutsch. Math.-Verein.\ \textbf{9} (1901), 115--121. Also in: Ges.\ Abh., vol.~II, pp.~122--127, Teubner, Leipzig, 1911.

\bibitem{Niven}
I.~Niven, \emph{Maxima and minima without calculus}, Dolciani Mathematical Expositions, no.~6, Mathematical Association of America, 1981.

\bibitem{MR94k:52003}
C.~Peri, J.~M.~Wills and A.~Zucco, \emph{On Blaschke's extension of Bonnesen's inequality}, Geom.\ Dedicata \textbf{48} (1993), 349--357.

\bibitem{MR18:760e}
C.~M. Petty, \emph{On the geometry of the {M}inkowski plane}, Riv. Mat. Univ.
  Parma \textbf{6} (1955), 269--292.

\bibitem{MR19:432a}
R.~R.~Phelps, \emph{Convex sets and nearest points}, Proc.\ Amer.\ Math.\ Soc.\
  \textbf{8} (1957), 790--797.

\bibitem{Radon}
J.~Radon, \emph{{\"U}ber eine besondere {A}rt ebener {K}urven}, Ber. Verh.
  S{\"a}chs.\ Ges.\ Wiss.\ Leipzig.\ Math.-Phys.\ Kl.\ \textbf{68} (1916), 23--28.

\bibitem{Riemann}
B.~Riemann, \emph{\"Uber die Hypothesen, welche der Geometrie zu Grunde liegen}, 
Abhandlungen der K\"oniglichen Gesellschaft der Wissenschaften zu G\"ottingen, 13, 1868.

\bibitem{MR47:5732}
J.~J.~Sch{\"a}ffer, \emph{The self-circumferences of polar convex disks}, Arch.\ Math.\
  (Basel) \textbf{24} (1973), 87--90.

\bibitem{MR98m:46018}
P.~Sch{\"o}pf, \emph{Orthogonality and proportional norms}, Anz. \"Osterreich.
  Akad. Wiss. Math.-Natur. Kl. \textbf{133} (1996), 11--16 (1997).

\bibitem{SzNagy}
B.~Sz.~Nagy, \emph{\"Uber ein geometrisches Extremalproblem}, Acta Univ.\ Szeged, Acta Sci.\ Math.\ \textbf{9} (1940), 253--257.

\bibitem{MR23:A4052}
L.~Tam{\'a}ssy, \emph{Ein {P}roblem der zweidimensionalen {M}inkowskischen
  {G}eometrie}, Ann.\ Polon.\ Math.\ \textbf{9} (1960/1961), 39--48.


\bibitem{MR52:4138}
A.~C.~Thompson, \emph{An equiperimetric property of {M}inkowski circles}, Bull.\
  London Math.\ Soc.\ \textbf{7} (1975), 271--272.

\bibitem{MR97f:52001}
A. C. Thompson: \emph{Minkowski Geometry.}
Cambridge University Press, Cambridge 1996.

\bibitem{Wallen}
L.~J.~Wallen, \emph{Kepler, the taxicab metric, and beyond: An isoperimetric primer}, College Math.\ J. \textbf{26} (1995), 178--190.


\end{thebibliography}
\end{document}